\documentclass[12pt, draft]{amsart}
\usepackage{epsfig,amsmath}

\begin{document}

\numberwithin{equation}{section}

\newtheorem{thm}[subsection]{Theorem}
\newtheorem{lem}[subsection]{Lemma}
\newtheorem{cor}[subsection]{Corollary}
\newtheorem{prop}[subsection]{Proposition}
\newtheorem{obs}[subsection]{Observation}
\newtheorem{claim}[subsection]{Claim}

\theoremstyle{definition}
\newtheorem{definition}[subsection]{Definition}
\newtheorem{rem}[subsection]{Remark}

\newcommand{\R}{\mathbf{R}}
\newcommand{\US}{\mathbf{S}}

\def\g{\gamma}
\def\eps{\varepsilon}
\def\z{\zeta}

\long\def\symbolfootnote[#1]#2{
	\begingroup
	\def\thefootnote{\fnsymbol{footnote}}\footnote[#1]{#2}
	\endgroup}
	

\catcode`@=11 \@mparswitchfalse  

\newcounter{mnotecount}[page]
\renewcommand{\themnotecount}{\arabic{mnotecount}}

\newcommand{\mnote}[1]
{\protect{\stepcounter{mnotecount}}$^{\mbox{\footnotesize  $
      \bullet$\themnotecount}}$
\marginpar{\null\hskip-9pt\raggedright\tiny\em
\themnotecount:\! \it #1} }


\parskip 5pt
\parindent 0pt
\baselineskip 15pt


\author{Bruce Solomon}


%

\title[Central cross sections make \dots quadric]{Central cross-sections make surfaces of revolution quadric}
\maketitle

\begin{center}
	\textit{Dedicated to Sue Swartz}
\end{center}


\section{Introduction}\label{sec:intro}

Quadric surfaces of revolution---ellipsoids, cones, paraboloids, cylinders, 
and hyperboloids---are the most basic
nontrivial surfaces in $\,\R^{3}\,$.
Euler noted \cite{boyer} that one can
describe them all, up to rigid motion, 
by appropriately choosing constants $\,a,b,\,$ and $\,c\,$ 
in the simple quadratic equation
	\begin{equation}\label{eqn:quad}
		x^{2}+y^{2}=az^{2}+ bz+c.
	\end{equation}
The left side of this equation
measures distance (squared) from the $z$-axis---the axis of 
revolution---while the right sets the vertical ``profile'' of the surface.


Archimedes knew that any plane tilted sufficiently far from the axis
of such a surface cuts it in an {ellipse}, and Fermat apparently realized the converse: elliptical cross-sections make a surface of revolution quadric \cite{coolidge}. 
One can verify Fermat's observation using 
high-school algebra, but here we prove a far stronger statement: 
\emph{Cross-sections which are merely \emph{central} (in the sense we shortly explain) make a surface of revolution quadric.}

This fact has generalizations and consequences, some of which we
sketch in our concluding remarks,
but we leave the full treatment of such extensions 
to a forthcoming article. Here we focus on the claim italicized above, and proceed now to define some terms we need to state and prove it.

A $\,C^{1}\,$ \textbf{loop} in $\,\R^{3}\,$ is the image of a \emph{periodic} $\,C^{1}\,$
mapping $\,\R\to\R^{3}\,$ whose derivative never vanishes. 

We call a loop in $\,\R^{3}\,$ \textbf{central} if it has symmetry
with respect to reflection through some point. We call that point its 
\textbf{center} and say that the loop is \textbf{centered} there. The center of a loop
clearly coincides with the center of mass, so it is unique. Circles, ellipses, and many other loops have central symmetry, though of course most loops do not. For instance, no loop sufficiently close to a triangle is central. 


By a \textbf{$\,C^{1}\,$ surface of revolution}, we 
mean a surface $\,S\,$ definable, after a rigid motion, 
as the locus
\begin{equation}\label{eqn:SoR}
		x^{2}+y^{2}=F(z),\quad |z|<q.
\end{equation}
Here $\,F\,$ is
a strictly positive, differentiable ``profile'' function. 
To indicate that a surface of revolution has this 
description for some $\,q>0\,$ (which bounds the vertical extent
of the suface) and some $\,F:(-q,q)\to(0,\infty)\,$, 
we say the surface lies in \textbf{standard position}. 
Note that we allow $\,F(z)\to\infty\,$ 
as $\,|z|\to q\,$. This can cause technical complications
that we will avoid, when necessary, by focusing on the restricted surface
\begin{equation*}\label{eqn:sdelta}
	S_{\delta}:=\left\{(x,y,z)\in S\colon |z|<q-\delta\right\}.
\end{equation*}
We deem a plane $\,P\,$ \textbf{transverse} to $\,S\,$ if $\,P\cap S\ne\emptyset\,$, {and} $\,P\,$ never coincides with the tangent plane 
of $\,S\,$ at any point of their intersection. 

All planes of interest here take the ``graphical'' form
\begin{equation}\label{eqn:nhp}
	z=m_{1}x+m_{2}y+\beta.
\end{equation}
We call $\,m:=\sqrt{m_{1}^{2}+m_{2}^{2}}\,$ the \textbf{slope}
and $\,\beta\,$ the \textbf{intercept} of this plane, which
we henceforth denote by $\,{P_{m,\beta}}\,$. 
Actually, an entire circle of planes have 
slope $\,m\,$ and intercept $\,\beta\,$, but the rotational symmetry of our problem makes them equivalent for our purposes;
we may safely ignore the ambiguity.

\begin{thm}[Main result]\label{thm:main}
Suppose we have a $\,C^{1}\,$ surface of revolution $\,S\,$ in standard position, and for some $\,\mu>0\,$, every plane of slope $\,m<\mu\,$ that cuts $\,S\,$ in a loop does so 
in a \emph{central} loop. Then $\,S\,$ is quadric.
\end{thm}

\section{The argument}

We prepare to prove Theorem \ref{thm:main} with three simple lemmas.
The first merely records some basic facts about surfaces of revolution:

\begin{lem}\label{lem:basics}
Suppose a $\,C^{1}\,$ surface of revolution $\,S\,$ lies in standard position, and $\,\delta>0\,$. Then there exists $\,\mu>0\,$ 
such that $\,P_{m,\beta}\cap S_{\delta}\,$ is a $\,C^{1}\,$ loop lying in the slab $\,|z-\beta|<\delta\,$ whenever $\,m<\mu\,$ and
$\,|\beta|<q-2\delta\,$.
\end{lem} 

One proves this using the implicit function theorem, 
together with the boundedness of $\,F\,$ and $\,F'\,$ when $\,|z|<q-\delta\,$. The result seems geometrically evident, so 
we leave further details to the reader.

The next lemma puts an amusing gloss 
on the classical mean value theorem for a 
differentiable function $\,f\,$. 
The latter equates the slope of the chord between  
$\,(a,f(a))\,$ and $\,(b,f(b))\,$ to the derivative 
$\,f'(c)\,$ for \emph{some} 
$\,c\,$ between $\,a\,$ and $\,b\,$.  
If we insist that $\,c\,$ coincide always with the 
\emph{midpoint} of $\,a\,$ and $\,b\,$, 
it turns out that we make $\,f\,$ \emph{quadratic}: 

\begin{lem}\label{lem:1}
	Suppose $\,f\,$ is differentiable on an open interval 
	$\,I$, and for some $\,\eps>0\,$ and all $\,\zeta\in I$,
	we have
	\begin{equation}\label{eqn:dq}
		 f'(\zeta)={f(\zeta+t)-f(\zeta-t)\over 2t}
	\end{equation}
	whenever $\,|t|<\eps\,$ and $\,\zeta\pm t\in I\,$.
	Then $\,f\,$ is quadratic on $\,I$.
\end{lem}

\begin{proof}
	Our assumptions make
	the right-hand side of (\ref{eqn:dq})---and hence the
	left too---differentiable with respect to 
	$\,\zeta\,$ when $\,|t|<\eps\,$ and $\,\zeta\pm t\in I\,$. 	
	For such $\,t\,$, taking $\,d/d\zeta\,$ has the 
	effect of replacing$\,f\,$ by $\,f'\,$ throughout (\ref{eqn:dq}).
	It follows that $\,f''\,$, and by iteration, 
	\emph{all} derivatives of $\,f\,$, exist as 
	continuous functions on $\, I \,$.  
	
	Now multiply (\ref{eqn:dq}) 
	by $\,2t\,$ and, under the same harmless 
	restrictions on $\,t\,$,
	differentiate thrice, this time with respect to $\,t\,$.
	One gets
	\[
		0=f'''(\zeta+t)+f'''(\zeta-t)\ .
	\]
	Setting $\,t=0\,$ now shows that $\,f'''\equiv 0\,$ 
	on $\, I \,$.
\end{proof}

The conclusion of the next lemma should now look promising. 
Recall that $\,q>0\,$ measures the vertical extent of a surface of
revolution $\,S\,$ in standard position.

\begin{lem}\label{lem:2}
	Suppose $\,S\,$ is a $\,C^{1}\,$ surface of revolution in standard position,
	$\,m>0\,$, and that for some $\,\beta\in(-q,q)\,$ the intersection 
	$\,P_{m,\beta}\cap S\,$ is a \emph{central} loop centered at 
	height $\,\zeta\,$. Then  
	\[
		F'(\zeta) = {F(\zeta+t)-F(\zeta-t)\over 2t}
	\]
	whenever 
	$\,|t|<\sup\{z-\zeta\colon (x,y,z)\in P_{m,\beta}\cap S\}\,$.
\end{lem}

\begin{proof}
	Define $\,b:=-\beta\,$. The rotational symmetry of $\,S\,$ then 
	lets us assume our plane	$\,P_{m,\beta}\,$ takes the form
	\begin{equation}\label{eqn:P}
		z =m\,x - b\ ,
	\end{equation}
	and by using this equation to eliminate $\,x\,$ in (\ref{eqn:SoR}),
	we can characterize our loop $\,P_{m,\beta}\cap S \,$,
	in the $(y,z)$-coordinate system on $\,P_{m,\beta}\,$, 
	as the locus
	\begin{equation}\label{eqn:gam}
	\left({b+z\over  m}\right)^{2}+ y^{2}=F(z)\,.
	\end{equation}
	Clearly, the reflection
	\begin{equation}\label{eqn:ySym}
		(y,z) \longmapsto\left(-y,\ z\right)\  
	\end{equation}
	preserves this loop, and thus $\,y=0\,$ at its center. 
	It follows that whenever a point with $\,(y,z)\,$ coordinates $\,(y,\,\z+t)\,$ satisfies
	(\ref{eqn:gam}), the point $\,(-y,\,\z-t)\,$ does too, so that  $\,y\,$ and $\,t\,$
	satisfy the simultaneous equations
	\begin{eqnarray*}\label{eqn:simul}
		F(\z+t) &=&\left(\frac{\bar b+ t}{m}\right)^{2}+y^{2}\\
		F(\z-t) &=&\left(\frac{\bar b- t}{m}\right)^{2}+(-y)^{2}\,,\\
	\end{eqnarray*}
	where $\,\bar b := b+\z\,$.
	Subtract the second equation from the first, simplify, and divide by
	$\,2t\,$ to get
	\begin{equation*}
		{F(\zeta+t)-F(\zeta-t)\over 2t}=\frac{2\bar b}{\,m^{2}}\,.
	\end{equation*}
	Letting $\,t\to 0\,$, 
	we see that the constant on the right must
	equal $\,F'(\zeta)\,$, and this proves the lemma.
\end{proof}

We can now verify our main result.

\begin{proof}[\emph{\textbf{Proof of Theorem \ref{thm:main}}}] 
	
	We must show that $\,F\,$ is quadratic. To do so, 	
	let $\,\delta>0\,$. Then our assumptions,
	together with Lemma \ref{lem:basics},
	guarantee that for some small but positive slope $\,m>0\,$,
	the plane $\,P_{m,\beta}\,$ given by
	\[
		z=mx+\beta
	\]
	cuts $\,S_{\delta}\,$ in a central
	loop lying in the slab $\,|z-\beta|<\delta\,$,
	provided only that $\,|\beta|<q-2\delta\,$. 
	So if we define a function
	\[
		\zeta:(-q+2\delta,\,q-2\delta)\to \R
	\]
	by making $\,\zeta(\beta)\,$ equal 
	the height of the center of the loop 
	$\,P_{m,\beta}\cap S_{\delta}\,$,
	the image of this function must contain the
	entire interval $\,|z|<q-3\delta\,$.
	
	Note also that the continuity and positivity
	of $\,F\,$ for $\,|z|<q\,$ ensure that
	\[
		\phi(\delta)
		:=\inf\left\{\sqrt{F(z)}\colon |z|<q-\delta\right\}>0\ .
	\]
	This means in particular that when $\,|\beta|<q-2\delta\,$, 
	the loop $\,P_{m,\beta}\cap S_{\delta}\,$ lies outside the 
	cylinder $\,x^{2}+y^{2}=\phi(\delta)^{2}\,$, so that the
	extreme values of $\,z\,$ on this loop  
	differ by at least $\,2m\,{\phi(\delta)}\,$. These extrema
	then deviate from $\,\zeta(\beta)\,$---the height 
	of the {center}---by at least
	\[ 
		\eps:=m\,{\phi(\delta)}\ .
	\]
	Lemma \ref{lem:2} now ensures that on the interval 
	$\,|\beta|<q-3\delta\,$,
	our profile function $\,F\,$ satisfies the assumptions of 
	Lemma \ref{lem:1} with the value of $\,\eps\,$ just determined,
	making $\,F\,$ quadratic on this interval.	
	But $\,\delta\,$ was arbitrary, so $\,F\,$ is 
	quadratic for all $\,|z|<q\,$, 
	and hence $\,S\,$ is quadric, as claimed.
\end{proof}  

\begin{rem}
	Though stated in $\,\R^{3}\,$, 
	all arguments above, and indeed our main theorem,
	generalize immediately to higher dimensions.
	One simply introduces coordinates
	\[
		(x_{1},\dots, x_{n-2},y,z)\in\R^{n}\ ,
	\]
	and replaces $\,x^{2}\,$ by 
	$\,|x|^{2}:=\sum_{i=1}^{n-2}x_{i}^{2}\,$ in (\ref{eqn:SoR})
	to define \textbf{hypersurface} \emph{of revolution in standard 
	position}. Everything above then generalizes to $\,\R^{n}\,$,
	with the word ``surface'' replaced by ``hypersurface'' throughout.
\end{rem}

\begin{rem}\label{rem:tubes}
	In a forthcoming paper, we apply the main result here in an 
	essential way to prove a far more general result. Roughly 
	speaking, we show there that any ``tube'' in 
	$\,\R^{3}\,$ that has compact, convex planar cross-sections, 
	all of them central, must be either quadric, or a cylinder
	over a centrally symmetric plane loop. 
\end{rem}
\goodbreak

\begin{rem}\label{rem:skew}
	A \textbf{skewloop} is a $\,C^{1}\,$ loop in $\,\R^{3}\,$ with no
	pair of parallel tangent lines. 
	The term was coined in \cite{gs}, which goes on to show 
	that convex quadrics are the only \emph{positively} curved 
	surfaces without skewloops. Because our main theorem says that 
	every \emph{non}-quadric surface of revolution has a
	\emph{non}-central cross-section, one can exploit the 
	``grafting'' technique from \cite[\S5]{gs} to show:
		\begin{quotation}\textit{	\!\!\!Every non-quadric surface 
			of revolution in $\,\R^{3}\,$ admits a skewloop.}
		\end{quotation}
	Conversely, by \cite{ss} (or, generically, \cite{tab}),
	no quadric admits a skewloop.
	So the main theorem here yields a characterization of the 
	one-sheeted hyperboloid 
	(the case where $a>0\,$ in (\ref{eqn:quad})) 
	as the only \emph{negatively} curved surface of revolution
	without skewloops. The generalization described in Remark 
	\ref{rem:tubes} will remove the need to assume rotational
	symmetry. Still, the main theorem
	here yields a first negatively curved counterpart 
	to the characterization of positively curved quadrics in \cite{gs}. 
\end{rem}


\section*{Acknowledgments} We thank the Lady Davis Foundation, the Technion (Israel Institute of Technology) and Indiana University for supporting this work.

\bigskip

\textbf{Bruce Solomon} 
began learning geometry from Barrett O'Neill while an undergraduate 
at UCLA. He earned a Ph.D. in the subject under F.~J. Almgren at
Princeton in 1982, and has practiced it at Indiana University,
Bloomington since 1983.

\textit{Math Department, Indiana University, 
Bloomington, IN 47405}\\
\textit{solomon@indiana.edu}
\end{document}